\documentclass{article}
\usepackage{amssymb,epsfig,graphicx}
\newcommand{\be}{\begin{equation}}
\newcommand{\ee}{\end{equation}}
\newcommand{\bea}{\begin{eqnarray}}
\newcommand{\eea}{\end{eqnarray}}
\newcommand{\barr}{\begin{array}}
\newcommand{\earr}{\end{array}}

\newcommand{\bpar}{\be \left\{ \begin{array}{lll}}
\newcommand{\epar}{ \end{array}\right. \ee }
\newcommand{\eparn}{ \end{array} \right.}

\newcommand\bvp{boundary value problem }
\newcommand\bvps{boundary value problems }

\newcommand\leftmat{\left(\begin{array}{cc}}
\newcommand\rightmat{\end{array}\right)}
\newcommand\leftvec{\left(\begin{array}{c}}
\newcommand\rightvec{\end{array}\right)}
\newcommand\E{{\rm e}^{ikL}}

\begin{document}
\bibliographystyle{plain}
\begin{flushright}
{\small to appear in  {\em IMA J. Appl. Math.}}
 \end{flushright}

\vspace{3mm}
\begin{center}
\Large
{\bf A transform method for linear evolution PDEs on a finite interval}

\vspace{3mm}

\large
A.S. Fokas$^{*}$ and B. Pelloni$^{**}$

\vspace{5mm}

*Department of Applied Mathematics and Theoretical Physics

Cambridge University

Cambridge CB3 0WA, UK.

{\em t.fokas@damtp.cam.ac.uk}

\vspace{2mm}

**Department of Mathematics

University of Reading

Reading RG6 6AX, UK

{\em b.pelloni@rdg.ac.uk}

\vspace{3mm}
\today
\end{center}
\normalsize

\begin{abstract}
We study initial boundary value problems for linear scalar evolution partial
differential equations, with spatial
derivatives of {\em arbitrary order}, posed on the domain
$\{t>0,\;0<x<L\}$.
We  show that the solution can be expressed
as an  integral in the
complex $k$-plane. This
integral is
defined in terms of
an $x$-transform of the initial condition and a  $t$-transform
of the boundary conditions. The derivation of this integral representation relies on the analysis of the {\em global relation},
which is an algebraic relation defined in the complex $k$-plane
coupling all boundary values of the solution.

For particular cases, such as the case of periodic boundary
conditions, or the case of boundary value problems for {\em
even}
order PDEs, it is possible to obtain directly 
from the global relation an alternative representation for the solution, 
in the form of an infinite series.
We stress however that
there exist initial boundary value problems for which the only
representation is an integral which
{\em cannot} be written as an infinite series. An
example of such a problem is provided by the linearised version of the
KdV equation. Similarly, in general the solution of odd-order
linear initial boundary value problems on a finite interval {\em
cannot} be expressed in terms of an infinite series.
\end{abstract}


\normalsize
\section{Introduction}
An evolution partial differential equation in one space dimension is 
characterised by its
{\em symbol}, which we denote by $\omega(k)$. This means that a
particular solution of the equation is given by
$$
E_{k}(x,t)={\rm e}^{ikx-\omega(k)t},\quad {\rm any}\;k\in\,\mathbb C.
$$
Physically significant examples of scalar evolution equations are:

\vspace{2mm}
{\bf (a)} the Schr\"odinger equation with zero potential
\be
iq_{t}+q_{xx}=0,\quad \omega(k)=ik^{2};
\label{eq1}
\ee
{\bf (b)} the heat equation
\be
q_{t}-q_{xx}=0,\quad \omega(k)=k^{2};
\label{eq2}
\ee
{\bf (c)} the Stokes equation
\be
q_{t}+q_{xxx}+q_{x}=0,\quad \omega(k)=i(k-k^{3}).
\label{eq3}
\ee

We note that equations (\ref{eq1}) and (\ref{eq3}) are the linearised
versions of the nonlinear Schr\"odinger and of the Korteweg-deVries
equations respectively.

\subsection*{ Notations}
The following notations, used throughout the paper, refer to a linear 
evolution PDE of order $n$, whose symbol $\omega(k)$ is assumed to be a 
polynomial 
of degree $n$.
\begin{itemize}
\item[(i)]
 $q_{0}(x)$ denotes the given initial condition, and
$\hat{q}_{0}(k)$ the Fourier transform of $q_{0}(x)$:
\be
q_{0}(x)=q(x,0),\;\; 0<x<L;
\;\;\;\hat{q}_{0}(k)=\int_{0}^{L}q_{0}(x){\rm
e}^{-ikx}dx,\;\;k\in{\mathbb C}.
\label{q0}
\ee
\item[(ii)]
$\{f_{j}(t)\}_{0}^{n-1}$ and $\{g_{j}(t)\}_{0}^{n-1}$ denote the
{\em boundary values} of the solution at $x=0$ and $x=L$ 
respectively, while
$\{\tilde{f}_{j}(t,k)\}_{0}^{n-1}$ and
$\{\tilde{g}_{j}(t,k)\}_{0}^{n-1}$ denote certain $t$-transforms of
$f_{j}$ and $g_{j}$:
\be
f_{j}(t)=\partial_{x}^{j}q(0,t),\;\;t>0; \;\;\;
\tilde{f}_{j}(t,k)=\int_{0}^{t}{\rm
e}^{\omega(k)s}f_{j}(s)ds,\;\;k\in{\mathbb C},\;t>0,
\label{fj}
\ee
\be
g_{j}(t)=\partial_{x}^{j}q(L,t),\;\;t>0; \;\;\;
\tilde{g}_{j}(t,k)=\int_{0}^{t}{\rm
e}^{\omega(k)s}g_{j}(s)ds,\;\;k\in{\mathbb C},\;t>0,
\label{gj}
\ee
where $j=0,\ldots,n-1$.
\item[(iii)]
The domain $D$ in the complex $k$-plane is defined by
\be
D=\{k\in{\mathbb C}:\;{\rm Re}\,\omega(k)\leq 0\}.
\label{Ddef}
\ee
$D^{+}$ and $D^{-}$ denote the part of $D$ in the upper and lower half
of the complex $k$-plane respectively:
\be
D^{+}=\{k\in{\mathbb C}:k\in\,D, \;{\rm Im}(k)\geq 0\},\quad
D^{-}=\{k\in{\mathbb C}:k\in\,D, \;{\rm Im}(k)\leq 0\}.
\label{D+-}
\ee
The oriented boundaries of $D^{+}$ and $D^{-}$ are denoted by $\partial
D^{+}$ and $\partial D^{-}$, where the orientation is such that the
interior of the domain $D$ is always on the left of the positive
direction.
\end{itemize}

\subsection*{Statement of the problem and assumptions}
{\em Let $q(x,t)$ satisfy a linear evolution equation with symbol
$\omega(k)$ in the domain
$$\left\{t>0,\quad 0<x<L\right\},$$
where $L$ is a finite positive constant. We assume that $\omega(k)$ 
is a
polynomial of degree $n$
such that the equation $\omega(k)=0$ has $n$ distinct roots, and
that Re\,$\omega(k)\geq 0$ for $k
\in\,{\mathbb R}$.
We assume that the initial
condition $q_{0}(x)$ is a given, sufficiently smooth, function.

We consider the two following questions:

{\bf (i)} Determine the number of boundary conditions that must be
prescribed at $x=0$ and $x=L$ in order to define a well posed problem.

{\bf (ii)} Given appropriate boundary conditions at the two ends of
the space interval, and assuming that these given functions have 
sufficient
smoothness and  are compatible with $q_{0}(x)$ at $x=0$ and
$x=L$, construct the solution $q(x,t)$.
}

\vspace{2mm} Problem (i) was solved in \cite{fsung2, pel}, where it was shown 
that  the number of boundary conditions that must be prescribed for a  well 
posed problem are $N$ at $x=0$ 
and $n-N$ at $x=L$, where
\be
N=\left\{\begin{array}{ll}
n/2& n\;\;{\rm even},\\
(n+1)/2 &n\;\;{\rm odd}, \;\;c_{n}>0,\\
(n-1)/2 &n\;\;{\rm odd}, \;\;c_{n}<0,
\end{array}
\right.
\label{Ndef}
\ee
where $c_{n}$ is the coefficient of $k^{n}$ in the symbol $\omega(k)$. 
We explain in appendix A  the motivation for this choice of  $N$.  

In this paper we address question (ii).

\subsection*{The new method}
The new method used here to analyse boundary value problems for
equations with spatial derivatives of {\em arbitrary} order $n$
is mathematically straightforward, yet it yields
results which are difficult to obtain by the
standard approaches. This method, which is the implementation to this
class of problems of the general approach introduced by
one of the authors \cite{fok}, involves the  steps outlined below.

\vspace{2mm}
{\it (a) Reformulation of the PDE}

A given evolution PDE with symbol  $\omega(k)$ can be written in
the form 
\be \left({\rm e}^{-ikx+\omega(k)t}q\right)_{t}-
\left({\rm e}^{-ikx+\omega(k)t} X\right)_{x}=0, \quad k\in{\mathbb
C}, \label{consform0} \ee 
where $q(x,t)$ is a solution of the PDE,
and the function $X(x,t,k)$ is given by the formula \be
X(x,t,k)=\sum_{j=0}^{n-1}c_{j}(k)\partial_{x}^{j}q(x,t),
\label{bigX} \ee where the coefficients $c_{j}(k)$ are known
polynomials in $k$. For example, for equation (\ref{eq1}), $n=2$,
$c_{1}(k)=-k$ and $c_{2}(k)=i$. The explicit form of $c_{j}(k)$
for an arbitrary $\omega(k)$ is given in section \ref{sec2}.

Equation (\ref{consform0}) is in the form $$
\frac{\partial P}{\partial t}-\frac{\partial Q}{\partial 
x}=0,
$$
with
$$
P(x,t)={\rm e}^{-ikx+\omega(k)t}q(x,t),\quad Q(x,t)={\rm 
e}^{-ikx+\omega(k)t}X(x,t,k).$$
Green's theorem applied 
to the domain ${\cal D}=\{[0,L]\times [0,t]\}$ yields
$$
\int\int_{\cal D}\left(\frac{\partial P}{\partial t}-\frac{\partial Q}{\partial 
x}\right)dx dt=\int_{\partial \cal D}Qdt+Pdx,
$$
thus
$$
\int_{\partial \cal D}Qdt+Pdx=0.
$$
Substituting in the latter expression the explicit definition of $P$ and 
$Q$, with $X(x,t,k)$ given by (\ref{bigX}),  we obtain the {\em global 
relation}
\be
\sum_{0}^{n-1} c_{j}(k)\left(\tilde{f}_{j}(t,k)-{\rm
e}^{-ikL}\tilde{g}_{j}(t,k)\right)=\hat{q}_{0}(k)-{\rm
e}^{\omega(k)t}\hat{q}(t,k),\quad k\in{\mathbb C},
\label{gr1}
\ee
where $\hat{q}_{0}$, $\tilde{f}_{j}$ and $\tilde{g}_{j}$ are defined
in the notations (equations (\ref{q0})-(\ref{gj}))and $\hat{q}(t,k)$
denotes the $x$-Fourier transform of $q(x,t)$.

Solving equation (\ref{gr1}) with respect to $\hat{q}(t,k)$ and then 
taking the inverse Fourier transform of the resulting expression, we obtain the 
following formula for $q(x,t)$:
\bea
q(x,t)&=&\frac{1}{2\pi}\left\{\int_{-\infty}^{\infty}{\rm
e}^{ikx-\omega(k)t}\hat{q}_{0}(k)dk-
\int_{-\infty}^{\infty}{\rm e}^{ikx-\omega(k)t}\sum_{j=0}^{n-1}
c_{j}(k)\tilde{f}_{j}(t,k)dk\right.
\nonumber \\
&+&\left.
\int_{-\infty}^{\infty}{\rm e}^{ik(x-L)-\omega(k)t}\sum_{j=0}^{n-1}
c_{j}(k)\tilde{g}(t,k)dk\right\}.
\label{intrepR}
\eea
We note however that this expression is {\em not 
effective}, since it contains the $t$-transforms of
{\em all} the boundary values of the solution $q(x,t)$, while 
only a subset of 
these boundary values is prescribed as boundary conditions.

\vspace{2mm}
{\it (b) The integral representation of the solution}

\begin{figure}
\begin{center}
\includegraphics[angle=360]{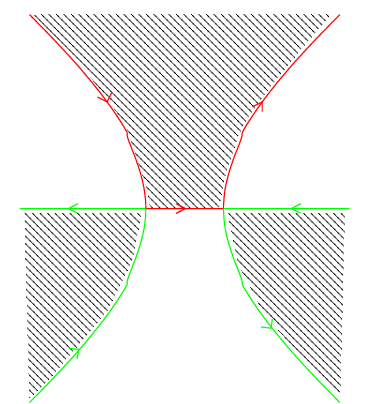}
\caption{The domain $D$ and the contours $\partial D^{+}$ (in red) 
and $\partial D^{-}$ (in green)  for
equation (\ref{eq3}). Note that $D^{-}$ has two connected components.}
\label{fig1}
\end{center}
\end{figure}

Using the analyticity properties of the functions $\tilde{f}_{j}$ and 
$\tilde{g}_{j}$ and Cauchy's theorem to deform the contour of 
integration, we show in section \ref{sec4} that
the expression (\ref{intrepR}) can be written as
\bea
q(x,t)&=&\frac{1}{2\pi}\left\{\int_{-\infty}^{\infty}{\rm
e}^{ikx-\omega(k)t}\hat{q}_{0}(k)dk-
\int_{\partial D^{+}}{\rm e}^{ikx-\omega(k)t}\sum_{j=0}^{n-1}
c_{j}(k)\tilde{f}_{j}(t,k)dk\right.
\nonumber \\
&-&\left.
\int_{\partial D^{-}}{\rm e}^{ik(x-L)-\omega(k)t}\sum_{j=0}^{n-1}
c_{j}(k)\tilde{g}(t,k)dk\right\},
\label{intrep1}
\eea
where  $\partial D^{+}$, $\partial D^{-}$ are defined in
the notations. The advantage of this form of the representation 
is that the $t$-transforms of all boundary values can now be 
computed explicitly. Indeed, it
will be shown in section \ref{sec4} that the functions $\tilde{f}_{j}$ and 
$\tilde{g}_{j}$, $j=0,\ldots,n-1$, 
for $k\in D^{\pm}$ respectively, can be expressed in terms of the 
given initial and boundary conditions.
For example, it will be shown in that section that the solution of
the \bvp\, for equation (\ref{eq3})  with the boundary conditions
\be
q(0,t)=f_{0}(t), \quad q(L,t)=g_{0}(t),\quad
q_{x}(L,t)=g_{1}(t),\quad t>0,
\label{bcc}
\ee
 is given by equation (\ref{intrep1}),
where the contours $\partial D^{+}$ and $\partial D^{-}$ are shown in
figure \ref{fig1}, $\omega(k)=i(k-k^{3})$, $c_{0}=k^{2}-1$, $c_{1}=-ik$, 
$c_{2}=1$, $\tilde{f}_{0}$, $\tilde{g}_{0}$ and $\tilde{g}_{1}$ are 
defined in terms of $f_{0}$, $g_{0}$ and $g_{1}$ by equations 
(\ref{fj}) and (\ref{gj}),
and the functions $\tilde{f}_{1}$, $\tilde{f}_2$ and
$\tilde{g}_{2}$ are given in terms of the given  initial and boundary
conditions by the following expressions:
\bea
i\tilde{f}_{1}&=&\frac{1}{\Delta(k)} \left[{\rm
e}^{-ikL}(N(\lambda_{2},t)-N(\lambda_{1},t))+{\rm
e}^{-i\lambda_{1}L}(N(k,t)-N(\lambda_{2},t))+{\rm
e}^{-i\lambda_{2}L}(N(\lambda_{1},t)-N(k,t))\right],
\nonumber \\
\tilde{f}_{2}&=&\frac{1}{\Delta(k)} \left[{\rm
e}^{-ikL}(\lambda_{2}N(\lambda_{1},t)-\lambda_{1}N(\lambda_{2},t))+{\rm
e}^{-i\lambda_{1}L}(kN(\lambda_{2},t)-\lambda_{2}N(k,t))
\right.\nonumber \\
&&\left.\hspace{30pt}+{\rm
e}^{-i\lambda_{2}L}(\lambda_{1}N(k,t)-kN(\lambda_{1},t))\right],
\nonumber \\
\tilde{g}_{2}&=&\frac{1}{\Delta(k)} \left[(N(k,t)(\lambda_2-\lambda_1)
+N(\lambda_{1},t)(\lambda_2-k)+N(\lambda_{2},t))(k-\lambda_1)\right],
\label{g23}
\eea
where
$\Delta(k)$ and $N(k)$ are defined by
\be
\Delta(k)={\rm e}^{-ikL}(\lambda_{1}-\lambda_{2})+{\rm
e}^{-i\lambda_{2}L}(k-\lambda_{1})+{\rm
e}^{-i\lambda_{1}L}(\lambda_{2}-k),
\label{detc1}
\ee
\be
N(k,t)=(1-k^{2})(\tilde{f}_{0}(t,k)-{\rm e}^{-ikL}\tilde{g}_{0}(t,k))-
ik{\rm e}^{-ikL}\tilde{g}_{1}(t,k)+\hat{q}_{0}(k),
\label{N3}\ee
and $\lambda_{1}(k)$, $\lambda_{2}(k)$ are the two roots of 
the polynomial $\lambda^{2}+\lambda k+k^{2}-1=0$.

\vspace{2mm}
{\it (c) The global relation and its analysis}

Although the  derivation of the global relation (\ref{gr1}) is elementary, 
this equation plays a central 
role in the analysis. The
crucial observation is that the functions $\tilde{f}_{j}$ and
$\tilde{g}_{j}$ depend on $k$ only through $\omega(k)$. 
Thus these functions are {\em invariant} under any transformation of
the complex $k$-plane
that leaves $\omega(k)$ invariant. These transformations are 
determined
by the roots of the equation
$
\omega(k)=\omega(\lambda)$. These roots are by our assumption 
distinct, and are given by
\be
\lambda_{0}(k)=k,
\lambda_{1},\ldots,\lambda_{n-1}.
\label{roots}
\ee
Replacing $k$ by $\lambda(k)$ in equation (\ref{gr1}) we obtain a
system of $n$ equations
\be
\sum_{j=0}^{n-1} c_{j}(\lambda_{l}(k))\left(\tilde{f}_{j}(t,k)-{\rm
e}^{-i\lambda_{l}(k)L}\tilde{g}_{j}(t,k)\right)=
\hat{q}_{0}(\lambda_{l}(k))-{\rm
e}^{\omega(k)t}\hat{q}(t,\lambda_{l}(k)),\quad
l=0,..,n-1.
\label{gr1n}
\ee
Ignoring for the moment the  {\em unknown} function
$\hat{q}(t,k)$, equations (\ref{gr1n}) can be considered as $n$
equations coupling the $2n$ unknown functions $\{\tilde{f}_{j},
\tilde{g}_{j}\}_{0}^{n-1}$;  $n$ of these functions can be computed 
immediately, and the remaining $n$ unknown functions can be obtained by solving
this system of $n$ equations. 

More specifically, let $f_{p}(t)$ and $g_{r}(t)$
denote the prescribed boundary conditions, where $p$ and $r$ take $N$
and $n-N$ integer values respectively. Then
$\tilde{f}_{p}$ and $\tilde{g}_{r}$ can be computed immediately (see
equations (\ref{fj}) and (\ref{gj})). Let $f_{P}(t)$ and $g_{R}(t)$
denote the remaining unknown boundary values. Solving the $n$ 
algebraic
equation (\ref{gr1n})  for $\tilde{f}_{P}$ and $\tilde{g}_{R}$ it
follows that these unknown functions can be expressed in terms of
$\hat{q}_{0}$, $\{\tilde{f}_{p},\tilde{g}_{r}\}$, and of a term
involving $\hat{q}(t,\lambda(k))$ and $1/\Delta(k)$, where
$\Delta(k)$ denotes the determinant of the relevant system. 
It was shown in \cite{fsung2, pel} that in order for the terms 
involving the unknown functions $\hat{q}(t,\lambda_{j}(k))$ {\em not} 
to contribute to $q(x,t)$, $N$ must be chosen by equation (\ref{Ndef}). 
Indeed, in this case, 
the terms involving
$\hat{q}(t,\lambda(k))$ appearing in the representation
of $\tilde{f}_{P}$ are bounded as $k\to\infty$, $k\in\,D^{+}$.
Therefore, if $\Delta(k)\neq 0$, $k\in\,D^{+}$,
using Cauchy theorem in the domain $D^{+}$, it follows that these 
expressions
give a zero contribution.
If $\Delta(k)=0$, for infinitely many $k\in\,D^{+}$, it can be shown 
that these zeros must be on $\partial D^+$ \cite{pel2}.
In this case, the contour deformation from $\mathbb R$ to $\partial 
D^{+}$ must avoid these zeros. The integral representation is then 
obtained along the indented contour, and as in the previous case the 
terms involving $\hat{q}(t,\lambda(k))$
give a zero contribution.

The analysis of
the terms $\tilde{g}_{R}$ is similar.

Thus even if $\Delta(k)$  has zeros in $D$, these zeros can be 
avoided by a contour indentation, and an integral representation of 
$q(x,t)$ can always be constructed. The determination of the zeros 
of $\Delta(k)$
%
becomes increasingly more complicated as the order $n$ of the 
equation grows. However, $\Delta(k)$  is always an analytic 
function of finite order, in the form of a finite exponential sum (the {\em order} of an entire
 function is a measure of its rate of growth as $k\to\infty$).  For such entire
 functions there exists an extensive theory, implying in particular 
 that such a function has
 infinitely many zeros accumulating at infinity, which lie along 
 specific rays in the complex plane, see \cite{levin}. A brief 
 overview of the relevant theory is given in Appendix B.
 This knowledge 
 is sufficient to determine when these zeros are inside or outside the 
 domain $D$.
 

\vspace{2mm}
{\it (d) The series representation}

For even order problems, separation of variables gives rise to a 
self-adjoint $x$-differential operator, which in this case can be 
used to construct the solution in the form of a series. For odd order 
problems, the associated $x$-operator is not, in general, 
self-adjoint. For {\em particular boundary conditions} (which do 
{\em not} include the uncoupled boundary conditions arising in many 
applications), it is possible to construct a self-adjoint extension, 
and thus for these particualr \bvps it is also possible to obtain the 
solution in the form of  a series.

There exist the following alternative, simple way to construct these 
series representations: {\em evaluate equations (\ref{gr1n}) at the zeros 
of $\Delta(k)$}. We note that this evaluation is possible only if all 
the terms appearing in equations (\ref{gr1n}) are bounded at these 
zeros. In this respect we differentiate between the even and odd order 
problems: for even order problems, this evaluation is always possible, 
independently of whether the zeros of $\Delta(k)$ are inside or 
outside of the domain $D$.   However, for odd order problems, it 
turns out that this evaluation is possible only if the zeros of 
$\Delta(k)$ are in $D$ (corresponding to  the case of odd order 
problems for 
which there exists a self-adjoitn extension). For typical \bvps, with 
uncoupled boundary conditions, the zeros of $\Delta(k)$ are {\em 
outside} $D$, and thus it is {\em not} possible to obtain a series 
representation (1.14).

Consider for example equation (\ref{eq1}) and assume for
simplicity that homogeneous Dirichlet boundary conditions are
prescribed. Then $\tilde{f}_{0}=\tilde{g}_{0}=0$ and
$\lambda_{1}=-k$, hence equations (\ref{gr1n}) become
\bea
i\tilde{f}_{1}-i{\rm e}^{-ikL} \tilde{g}_{1}&=&\hat{q}_{0}(k)-{\rm
e}^{ik^{2}t}\hat{q}(t,k),\nonumber
\\
i\tilde{f}_{1}-i{\rm e}^{ikL} \tilde{g}_{1}&=&\hat{q}_{0}(-k)-{\rm
e}^{ik^{2}t}\hat{q}(t,-k).
\nonumber
\eea
Subtracting these equations we obtain
$$
i({\rm e}^{ikL}-{\rm 
e}^{-ikL})\tilde{g}_{1}=\hat{q}_{0}(k)-\hat{q}_{0}(-k)
-{\rm
e}^{ik^{2}t}\left(\hat{q}(t,k)-\hat{q}(t,-k)\right).
$$
Evaluating this equation at the values of $k$ for which the
coefficient of $\tilde{g}_{1}$ vanishes, i.e. 
$k=k_{m}:=\frac{m\pi}{L}$,
$m\in{\mathbb Z}$, and using the definition of $\hat{q}_{0}(k)$ and
$\hat{q}(t,k)$ we find
\be
\hat{q}^{(sin)}(t,k_{m})={\rm
e}^{-ik_{m}^{2}t}\hat{q}_{0}^{(sin)}(k_{m}),\quad
\hat{q}^{(sin)}(k_{m})=\int_{0}^{L}\sin(k_{m}x)q(x)dx.
\label{sintr}
\ee
Equation (\ref{sintr}a) can be inverted by the well known formula and
then $q(x,t)$ is expressed in the form of a sine series.
Equation (\ref{sintr}) can also be obtained by the use of the $x$-sine 
transform, which is the appropriate $x$-transform for this problem.  

The above computation relies on  the fact that all 
terms involved in the global relation are bounded at the zeros of 
$\Delta(k)$. This is to be contrasted with the case of
equation (\ref{eq3}) with $q(0,t)=q(L,t)=q_x(L,t)=0$. In this case, 
letting
$\tilde{f}_{0}=\tilde{g}_{0}=\tilde{g}_{1}=0$, and evaluating the 
global relation at $k$, $\lambda_{1}(k)$ and $\lambda_{2}(k)$, 
we obtain a system of three algebraic equations for 
$\tilde{f}_{1}$, $\tilde{f}_{2}$ and $\tilde{g}_{2}$. The determinant
of this system vanishes at $k=k_{m}$, $m\in\,{\mathbb Z}$, where for 
large $k$, arg$(k_{m})$ is either $\pi/6$, $5\pi/6$ or $3\pi/2$, 
so that $k_m \notin D$. It can be verified that 
the functions
$\tilde{f}_{1}$, $\tilde{f}_{2}$ and $\tilde{g}_{2}$ involve ${\rm 
e}^{i(k-k^{3})t}$, thus they become 
{\em 
unbounded} at $k=k_{m}$, $k\to \infty$. Hence, we {\em cannot} evaluate this system 
at the zeros of its determinant. A more detailed discussion is given 
in \cite{pel2}.

We note that when such a representation exists, it 
can also be obtained by using the explicit residue computation of the 
general integral representation.

\section{The elements of the method}\label{sec2}
\setcounter{equation}{0}
It can be verified (see \cite{ima1}) that a PDE with
symbol $\omega(k)$ can be written in the form
\be
\left({\rm e}^{-ikx+\omega(k)t}q\right)_{t}-
\left({\rm e}^{-ikx+\omega(k)t}
X\right)_{x}=0, \quad k\in{\mathbb C},
\label{consform}
\ee
where the function
$X(x,t,k)$ is given by the formula
\be
X(x,t,k)=i\frac{\omega(k)-
\omega(-i\partial_{x})}{k+i\partial_{x}}q(x,t).
\label{bigx}
\ee
In what follows we give the function $X$, the domains $D^{+}$, $D^{-}$,
and equations (\ref{gr1n})  for
equations (\ref{eq1})-(\ref{eq3}).

\subsubsection{(a) The equation (\ref{eq1})}
The symbol is
$$\omega(k)=ik^{2}.$$
Indeed, ${\rm
e}^{ikx-ik^{2}t}$ satisfies this equation. In this case 
$$
{\rm Re}\,\omega(k)={\rm Re}(ik^{2})=-2{\rm Re}(k){\rm Im}(k).
$$
Thus the domain $D$ is the union of the first and third 
quadrants of the $k$ complex plane:
\be
D=\{k\in{\mathbb C}^{+}:{\rm Re}(k){\rm Im}(k)\geq 0\}.
\label{D1}
\ee
The equation $\omega(k)=\omega(\lambda)$ implies 
$$
\lambda^{2}-k^{2}=(\lambda-k)(\lambda+k)=0$$ hence $\lambda_{1}(k)=-k$.
Equation (\ref{bigX}) yields
$$
X(x,t,k)=i\frac{ik^{2}+i\partial_{x}^{2}}
{k+i\partial_{x}}q=-\frac{(k+i\partial_{x})(k-i\partial_{x})}
{k+i\partial_{x}}q=iq_{x}-kq.
$$
The PDE (\ref{eq1})  is thus equivalent to the expression
\be
\left({\rm e}^{-ikx+ik^{2}t}q\right)_{t}-
\left({\rm e}^{-ikx+ik^{2}t}
(iq_{x}-kq)\right)_{x}=0, \quad k\in{\mathbb C}.
\label{eq1cs}
\ee
The global relation is
\be
i\tilde{f}_{1}(t,k)-k\tilde{f}_{0}(t,k)-{\rm
e}^{-ikL}\left(i\tilde{g}_{1}(t,k)-k\tilde{g}_{0}(t,k)\right)
=\hat{q}_{0}(k)-{\rm
e}^{ik^{2}t}\hat{q}(t,k).
\label{agr1}
\ee
Since $\lambda_{1}=-k$, we supplement this equation with the equation
\be
i\tilde{f}_{1}(t,k)+k\tilde{f}_{0}(t,k)-{\rm
e}^{ikL}\left(i\tilde{g}_{1}(t,k)+k\tilde{g}_{0}(t,k)\right)=
\hat{q}_{0}(-k)-{\rm
e}^{ik^{2}t}\hat{q}(t,-k).
\label{sys1}
\ee

\subsubsection{(b) The equation (\ref{eq2})}
\begin{figure}
\begin{center}
 \includegraphics[angle=360]{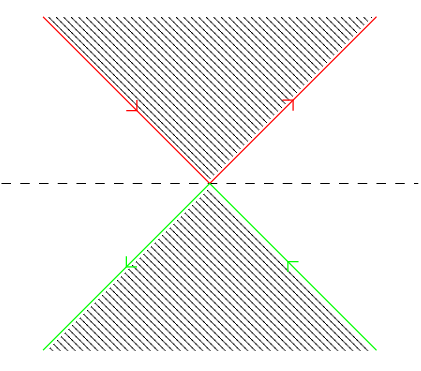}
\caption{The domain $D$
and the contours $\partial D^{+}$ (in red)  and $\partial D^{-}$ (in 
green)
for
equation (\ref{eq2})}
\label{heatdom}
\end{center}
\end{figure}

The symbol is
$$\omega(k)=k^{2}.$$
The domain $D$ is in this case given by 
\be
D=\{k\in{\mathbb C}:\;\pi/4\leq {\rm arg}\,(k)\leq 3\pi/4\;\;{\rm or}\;\;
5\pi/4\leq {\rm arg}\,(k)\leq 7\pi/4\}.
\label{hypb}
\ee
and $D^{\pm}=D\cap {\mathbb C}^{\pm}$, see figure \ref{heatdom}.

The equation $\omega(k)=\omega(\lambda)$ implies, as for example (a),
that  $\lambda_{1}(k)=-k$.
Equation (\ref{bigX}) yields
\be
X(x,t,k)=i\frac{k^{2}+\partial_{x}^{2}}
{k+i\partial_{x}}q=i\frac{(k+i\partial_{x})(k-i\partial_{x})}
{k+i\partial_{x}}q=ikq+q_{x}.
\label{bigX2}
\ee
The PDE (\ref{eq2}) is thus equivalent to the expression
\be
\left({\rm e}^{-ikx+k^{2}t}q\right)_{t}-
\left({\rm e}^{-ikx+k^{2}t}
(ikq+q_{x})\right)_{x}=0, \quad k\in{\mathbb C}.
\label{eq2cs}
\ee
The global relation is
\be
ik\tilde{f}_{0}(t,k)+\tilde{f}_{1}(t,k)-{\rm
e}^{-ikL}\left(ik\tilde{g}_{0}(t,k)+\tilde{g}_{1}(t,k)\right)
=\hat{q}_{0}(k)-{\rm
e}^{k^{2}t}\hat{q}(t,k).
\label{agr2}
\ee
Since $\lambda_{1}=-k$, we supplement this equation with the equation
\be
 -ik\tilde{f}_{0}(t,k)+\tilde{f}_{1}(t,k)-{\rm 
e}^{ikL}\left(-ik\tilde{g}_{0}(t,k)+\tilde{g}_{1}(t,k)\right)=
\hat{q}_{0}(-k)-{\rm
e}^{k^{2}t}\hat{q}(t,-k).
 \label{sys2}
 \ee

\subsubsection{(c) The equation (\ref{eq3})}

The symbol is
$$\omega(k)=i(k-k^{3}).$$
The domain $D$ is in this case
\be
D=\{k\in\,{\mathbb C}: \,{\rm Im}\,(k)\left[3{\rm Re}\,(k)^{2}-
{\rm Im}\,(k)^{2}-1\right]\leq 0\},
\label{hypc}
\ee
see figure \ref{fig1}.

The equation $\omega(k)=\omega(\lambda)$ implies
$$
(\lambda-\lambda^{3})-(k-k^{3})=(k-\lambda)(\lambda^{2}+\lambda 
k+k^{2}-1)=0,$$
hence
$$
\lambda_{1}=\frac{\sqrt{4-3k^{2}}}{2}-\frac{k}{2},\quad
\lambda_{2}=-\frac{\sqrt{4-3k^{2}}}{2}-\frac{k}{2}.$$

Equation (\ref{bigX}) yields
\bea
X(x,t,k)&=&i\frac{i(k-k^{3})-(\partial_{x}+\partial_{x}^{3})}
{k+i\partial_{x}}q=-\frac{(k+i\partial_{x}-k^{3}-i\partial_{x}^{3})}{k+i\partial_{x}}
q
\nonumber \\
&=&\left\{-1+\frac{(k+i\partial_{x})(k^{2}-\partial_{x}^{2}-ik\partial_{x})}
{k+i\partial_{x}}\right\}q=(k^{2}-1)q-ikq_{x}-q_{xx}.
\label{bigX3}
\eea
The PDE (\ref{eq3}) is thus equivalent to the expression
\be
\left({\rm e}^{-ikx+i(k-k^{3})t}q\right)_{t}-
\left({\rm e}^{-ikx+i(k-k^{3})t}
((k^{2}-1)q-ikq_{x}-q_{xx})\right)_{x}=0, \quad k\in{\mathbb C}.
\label{eq3cs}
\ee
The global relation is
\be
\left[(k^{2}-1)\tilde{f}_{0}-ik\tilde{f}_{1}-\tilde{f}_{2}\right]-{\rm
e}^{-ikL}\left[(k^{2}-1)\tilde{g}_{0}-ik\tilde{g}_{1}-\tilde{g}_{2}\right]
=\hat{q}_{0}(k)-{\rm
e}^{i(k-k^{3})t}\hat{q}(t,k).
\label{agr3}
\ee
Supplementing this equation with the equations obtained from it by 
replacing $k$ with $\lambda_{1}$ and 
$\lambda_{2}$,
we obtain 
\bea
\left[(\lambda_{1}^{2}-1)\tilde{f}_{0}-i\lambda_{1}\tilde{f}_{1
}-\tilde{f}_{2}\right]-{\rm
e}^{-i\lambda_{1}L}\left[(\lambda_{1}^{2}-1)\tilde{g}_{0}-
i\lambda_{1}\tilde{g}_{1}-\tilde{g}_{2}\right]
&=&\hat{q}_{0}(\lambda_{1})-{\rm
e}^{i(k-k^{3})t}\hat{q}(\lambda_{1},t),
\nonumber \\
\left[(\lambda_{2}^{2}-1)\tilde{f}_{0}-i\lambda_{2}\tilde{f}_{1}
-\tilde{f}_{2}\right]-{\rm
e}^{-i\lambda_{2}L}\left[(\lambda_{2}^{2}-1)\tilde{g}_{0}
-i\lambda_{2}\tilde{g}_{1}-\tilde{g}_{2}\right]
&=&\hat{q}_{0}(\lambda_{2})-{\rm
e}^{i(k-k^{3})t}\hat{q}(\lambda_{2},t).
\nonumber \\
\label{sys3}
\eea

%
%
%
%
%

\section{The zeros of $\Delta(k)$}\label{sec3}
We start by making the
following general observations, valid for an arbitrary $\omega(k)$:

{\bf (i)} The functions $\hat{q}_{0}(k)$ and $\hat{q}(t,k)$, which
are entire functions of $k$, are bounded as $k\to\infty$ in
${\mathbb C}^{-}$, the lower half of the complex $k$-plane.

{\bf (ii)} The functions ${\rm e}^{ikL}\hat{q}_{0}(k)$ and ${\rm
e}^{ikL}\hat{q}(t,k)$ are bounded as $k\to\infty$ in
${\mathbb C}^{+}$, the upper half of the complex $k$-plane.

{\bf (iii)}
The functions $\tilde{f}_{j}$ and $\tilde{g}_{j}$, which are entire
functions of $k$, are bounded as $k\to\infty$ if Re\,$\omega(k)\leq 
0$.

\subsubsection{Example (a)}
Since $n=2$ and $N=1$, one boundary condition must be prescribed at 
each end. 
We consider equations (\ref{agr1})-(\ref{sys1}) as two equations relating the four
terms $\tilde{f}_{0}$, $\tilde{f}_{1}$, $\tilde{g}_{0}$,
$\tilde{g}_1$, treating for the moment the function $\hat{q}(t,k)$
as a known function. If one of the functions $f_{j}$'s and one of the 
$g_{j}$'s are given, then
the determinant of the system is  of the 
form $\Delta(k)={\rm e}^{ikL}\pm{\rm e}^{-ikL}$. Hence its zeros are 
real, and given by either $k_{m}=m\pi/L$ or $k_{m}=(2m+1)\pi/2L$.
Evaluating the solution of system (\ref{agr1})-(\ref{sys1}) at these zeros 
yields the series representation of the solution for any 
such \bvp.

\vspace{2mm}
{\bf The Dirichlet problem}

The given boundary conditions are
\be
q(0,t)=f_{0}(t),\quad q(L,t)=g_{0}(t),\quad t>0.
\label{bca1}
\ee
Then $q(x,t)$ is given by 
\be
q(x,t)=\frac{i}{4L}\sum_{m}\sin(k_{m}x){\rm
e}^{-ik_{m}^{2}t}\left[N(k_m,t)-N(-k_m,t)\right],\quad
k_{m}=\frac{m\pi}{L},
\label{sin1}
\ee
where $N(k,t)$ is the known function
\be
N(k,t)=k\left(\tilde{f}_{0}(t,k) -{\rm e}^{-ikL}\tilde{g}_{0}(t,k)
\right)+\hat{q}_{0}(k).
\label{Nk}
\ee
Indeed, equations
(\ref{agr1})-(\ref{sys1}) yield
\bea
i\tilde{f}_{1}-{\rm e}^{-ikL}i\tilde{g}_{1}&=&N(k,t)-{\rm
e}^{ik^{2}t}\hat{q}(t,k),
\nonumber
\\
i\tilde{f}_{1}-{\rm e}^{ikL}i\tilde{g}_{1}&=&N(-k,t)-{\rm
e}^{ik^{2}t}\hat{q}(t,-k).
\label{A.2}
\eea
Subtracting equations (\ref{A.2}) we find
\be
({\rm e}^{ikL}-{\rm
e}^{-ikL})i\tilde{g}_{1}=N(k,t)-N(-k,t)
-{\rm
e}^{ik^{2}t}(\hat{q}(t,k)-\hat{q}(t,-k)),\quad k_{m}\in\,{\mathbb R}.
\label{lasteq}
\ee
We evaluate this equation at the values of $k$ for which ${\rm 
e}^{ikL}-{\rm
e}^{-ikL}=0$, i.e.
$$k=k_{m}=\frac{m\pi}{L}, \;m\in\,{\mathbb Z}.
$$
The definition of $\hat{q}(t,k)$ implies
$$
\int_{0}^{L}\sin(k_{m}x)q(x,t)dx=\frac{i}{2}{\rm
e}^{-ik_{m}^{2}t}\left[N(k_m,t)-N(-k_m,t)\right].
$$
Hence inverting this expression  we find (\ref{sin1}).

\subsubsection{Example (b)}

As for the previous example,  $n=2$ and $N=1$, one boundary condition must be prescribed at 
each end. 
The determinant 
of any boundary value problem obtained by prescribing one of the 
functions 
$f_{j}$'s and one of the $g_{j}$'s   is again of the 
form of example (a). The only 
difference is that for this 
example the real axis is outside the domain $D$ (except for the 
point $k=0$). 
However, still all terms appearing in (\ref{agr2})-(\ref{sys2}) are 
bounded for $k\in\,\mathbb R$, 
hence the same computation as before yields the series 
representation of the solution. 

\vspace{2mm}
{\bf The Dirichlet problem}

The given boundary conditions are
\be
q(0,t)=f_{0}(t),\quad q(L,t)=g_{0}(t),\quad t>0.
\label{bca2}
\ee
Then $q(x,t)$ is given by 
\be
q(x,t)=\frac{i}{4L}\sum_{m}\sin(k_{m}x){\rm
e}^{-k_{m}^{2}t}\left[N(k_m,t)-N(-k_m,t)\right],\quad
k_{m}=\frac{m\pi}{L},
\label{sin2}
\ee
where $N(k,t)$ is the known function
\be
N(k,t)=\hat{q}_{0}(k)-ik\left(\tilde{f}_{0}(t,k) -{\rm e}^{-ikL}\tilde{g}_{0}(t,k)
\right).
\label{Nk2}
\ee
The derivation is analogous to the case of example (a).

\vspace{2mm}
{\bf The Robin problem}

The given  boundary
conditions are now
\be
q_{x}(0,t)-\alpha q(0,t)=h_{1}(t),\quad
q_{x}(L,t)-\beta q(L,t)=h_{2}(t),\quad t>0,
\label{robin}
\ee
where $\alpha$ and $\beta$ are given real constants, $\alpha\neq 
\beta$.
%
Using the boundary conditions (\ref{robin}), the global
relation becomes
\be
(\alpha+ik)\tilde{f}_{0}-{\rm
e}^{-ikL}(\beta+ik)\tilde{g}_{0}=N(k,t)-{\rm
e}^{k^{2}t}\hat{q}(t,k),\quad k\in{\mathbb C},
\label{grheat2}
\ee
where the known function
$N(k,t)$ is given by
\be
N(k,t)=\hat{q}_{0}(k)+\int_{0}^{t}{\rm e}^{k^{2}s}\left({\rm
e}^{-ikL}h_{2}(s)-h_{1}(s)\right)ds.
\label{Nheat}
\ee
Since $\lambda_{1}=-k$, replacing $k$
by $-k$ in equation (\ref{grheat2}) we find
\be
(\alpha-ik)\tilde{f}_{0}-
{\rm e}^{ikL}(\beta-ik)\tilde{g}_{0}=N(-k,t)-{\rm
e}^{k^{2}t}\hat{q}(t,-k),\quad k\in{\mathbb C}.
\label{grheat3}
\ee
Solving the system (\ref{grheat2})-(\ref{grheat3}) we obtain
\bea
\Delta(k)\tilde{f}_{0}&=&{\rm e}^{ikL}(\beta-ik)N(k,t)-{\rm e}^{-ikL}(\beta+ik)N(-k,t)
\nonumber \\
&&-{\rm
e}^{k^{2}t}\left[{\rm e}^{ikL}(\beta-ik)\hat{q}(t,k)-{\rm e}^{-ikL}(\beta+ik)\hat{q}(t,-k)\right],
\nonumber \\
\label{solveheat1}
\\
\Delta(k)\tilde{g}_{0}&=&(\alpha-ik)N(k,t)-(\alpha+ik)N(-k,t)
-{\rm
e}^{k^{2}t}\left[(\alpha-ik)\hat{q}(t,k)-(\alpha+ik)\hat{q}(t,-k)\right],
\nonumber \\
\label{solveheat2}
\eea
where
\be
\Delta(k)=(\alpha+ik)(\beta-ik)\E-(\alpha-ik)(\beta+ik){\rm
e}^{-ikL}.
\label{heatdet}
\ee
Let $k_{m}$ be defined by
\be
k_{m}:\quad (\alpha+ik_{m})(\beta-ik_{m}){\rm
e}^{ik_{m}L}=(\alpha-ik_{m})(\beta+ik_{m}){\rm
e}^{-ik_{m}L}.
\label{zerosheat}
\ee
Evaluating equation (\ref{solveheat1})
at $k=k_{m}$, where
$\Delta(k_{m})=0$, and using
the definition of
$\hat{q}(t,k_{m})$, we find
$$
\int_{0}^{L}
\left[{\rm e}^{ikL}(\beta-ik_{m}){\rm e}^{-ik_{m}x}-{\rm e}^{-ikL}(\beta+ik_{m}){\rm
e}^{ik_{m}x}\right]q(x,t)dx=
$$
\be {\rm e}^{-k_{m}^{2}t}\left[
{\rm e}^{ikL}(\beta-ik_{m})N(k_{m},t)-{\rm e}^{-ikL}(\beta+ik_{m})N(-k_{m},t)\right].
\label{hfin} \ee 
Note that the function
 $\Delta(k)$ defined by
 (\ref{heatdet}) is entire and of finite order, hence it has infinitely
 many zeros accumulating at infinity. Moreover,
 this function has infinitely many zeros on the real axis, and
 all zeros are asymptotically on ${\mathbb R}$ \cite{pel2}. It follows that 
 the only one of these zero in $D$ is $k=0$, where the numerator also
 vanishes.
 However, equation (\ref{zerosheat}) cannot be solved explicitly 
for the $k_{m}$'s. Thus although the general theory implies that 
equation (\ref{hfin}) can be
solved for $q(x,t)$, it does not appear to yield an effective
representation for the solution of this problem. On the other
hand, the integral representation  {\em does}
provide an effective representation, see section \ref{sec4}.

\subsubsection{Example (c)}

In this case $n=3$ and $N=1$, thus for a  well posed problem 
for equation (\ref{eq3}) 
one boundary condition must be
prescribed
at $x=0$ and two boundary conditions must be prescribed at $x=L$. 

We consider equations (\ref{agr3})-(\ref{sys3}) as three equations relating the six
terms $\tilde{f}_{0}$, $\tilde{f}_{1}$, $\tilde{f}_{2}$, $\tilde{g}_{0}$,
$\tilde{g}_1$ and $\tilde{g}_{2}$, temporarily treating the function $\hat{q}(t,k)$
as a known function. If one of the functions $f_{j}$'s and two of the 
$g_{j}$'s are given, then
the determinant of the system is  always of the 
form
\be
\Delta(k)=K_{0}(k,\lambda_{1},\lambda_{2}){\rm e}^{-ikL}+
K_{1}(k,\lambda_{1},\lambda_{2}){\rm e}^{-i\lambda_{1}L}+
K_{2}(k,\lambda_{1},\lambda_{2}){\rm e}^{-i\lambda_{2}L}
\ee
where $K_{i}(k,\lambda_{1},\lambda_{2})$ are at most quadratic 
functions of the three arguments. In the limit as $k\to\infty$, up to 
multiple of $k$, this 
function behaves like the function
\be
\tilde{\Delta}(k)=\tilde{K}_{0}(\zeta){\rm e}^{-ikL}+
\tilde{K}_{1}(\zeta){\rm e}^{-i\zeta k L}+
\tilde{K}_{2}(\zeta){\rm e}^{-i\zeta^{2}kL},\quad \zeta={\rm 
e}^{2\pi i/3}.
\label{det3rd}
\ee
The particular form of the coefficients $\tilde{K}_{j}$ depends on the 
particular boundary conditions prescribed. For example if the 
prescribed conditions are given by (\ref{bcc}),
then 
\be
\tilde{\Delta}(k)=k\zeta(\zeta-1)\left[{\rm e}^{-ikL}+
\zeta {\rm e}^{-i\zeta k L}+\zeta^{2}{\rm e}^{-i\zeta^{2}kL}\right].
\label{tildedet3}
\ee
By the general theory presented 
in \cite{levin}, and briefly 
summarised in Appendix B, the infinitely many zeros of (\ref{det3rd}) 
depend only on the three exponentials appearing in (\ref{tildedet3}). 
These zeros
accumulate 
at infinity along the three lines arg$(k)=\pi/6$, arg$(k)=5\pi/6$ 
and arg$(k)=3\pi/2$, which are all outside $D$. 

We now show that in this case it is not possible to derive a 
series representation for the solution.
Equations  (\ref{agr3})-(\ref{sys3}) yield
\bea
-ik\tilde{f}_{1}-\tilde{f}_{2}+
{\rm 
e}^{-ikL}\tilde{g}_{2}
&=&N(k,t)-{\rm
e}^{i(k-k^{3})t}\hat{q}(t,k),
\nonumber \\
-i\lambda_{1}\tilde{f}_{1
}-\tilde{f}_{2}+{\rm
e}^{-i\lambda_{1}L}\tilde{g}_{2}
&=&N(\lambda_{1,t})-{\rm
e}^{i(k-k^{3})t}\hat{q}(\lambda_{1},t),
\nonumber \\
-i\lambda_{2}\tilde{f}_{1}
-\tilde{f}_{2}+{\rm
e}^{-i\lambda_{2}L}\tilde{g}_{2}
&=&N(\lambda_{2},t)-{\rm
e}^{i(k-k^{3})t}\hat{q}(\lambda_{2},t),
\nonumber \\
\label{sys3dir}
\eea
where
$N(k,t)$ is given by (\ref{N3}).
Solving the above system with respect e.g. to $\tilde{g}_{2}(t,k)$, we 
obtain
$$
\Delta(k)\tilde{g}_{2}(t,k)=\left[N(\lambda_{1},t)(\lambda_{2}-k)+
N(\lambda_{2},t)(k-\lambda_{1})+
N(k,t)(\lambda_{2}-\lambda_{1})\right]
$$
$$
-{\rm
e}^{i(k-k^{3})t}\left[\hat{q}(t,\lambda_{1})(k-\lambda_{1})+
\hat{q}(t,\lambda_{2})(\lambda_{2}-k)+
\hat{q}(t,k)(\lambda_{2}-\lambda_{1})\right].
$$
We cannot evaluate this expression at the zeros of 
$\Delta(k)$ as $k\to\infty$. For example, at the zeros which lie 
in $\mathbb C^-$ (i.e.
the zeros which asymptotically have argument equal to $3\pi/2$), 
the terms  $\hat{q}(t,\lambda_{1})$ and $\hat{q}(t,\lambda_{2})$
are {\em not bounded} as $k\to\infty$. 

\section{The integral representation of the solution}\label{sec4}
\setcounter{equation}{0}
We first derive equation (\ref{intrep1}). The global relation
(\ref{gr1}) yields
$$
\hat{q}(t,k)={\rm e}^{-\omega(k)t}\hat{q}_{0}(k)-{\rm
e}^{-\omega(k)t}\left[\tilde{f}(t,k)-{\rm 
e}^{-ikL}\tilde{g}(t,k)\right],$$
where
$$\tilde{f}(t,k)=\sum_{0}^{n-1}c_{j}(k)\tilde{f}_{j}(t,k),\quad
\tilde{g}(t,k)=\sum_{0}^{n-1}c_{j}(k)\tilde{g}_{j}(t,k),
$$
$$
\hat{q}(t,k)=\int_{0}^{L}{\rm e}^{-ikx}q(x,t)dx,\quad
\hat{q}_{0}(k)=\hat{q}(0,k).$$
Taking the inverse Fourier transform of $\hat{q}(t,k)$, we obtain
$$
q(x,t)=\frac{1}{2\pi}\int_{-\infty}^{\infty}{\rm
e}^{ikx-\omega(k)t}\left\{\hat{q}_{0}(k)-\tilde{f}(t,k)+{\rm
e}^{-ikL}\tilde{g}(t,k)\right\}dk.
$$
It follows from the  definition of $D$  that for
$k\notin D$ and for any $t>0$, the functions
${\rm e}^{-\omega(k)t}\tilde{f}(t,k)$ and
${\rm e}^{-\omega(k)t}\tilde{g}(t,k)$ are  bounded
as $k\to\infty$.
Thus

\vspace{2mm}
\begin{tabular}{ll}
${\rm e}^{ikx-\omega(k)t}\tilde{f}(k)$&is
analytic and bounded for $k\in\,{\mathbb C}^{+}\setminus D$,
 \\
 \\
 ${\rm e}^{ik(x-L)-\omega(k)t}\tilde{g}(k)$&is
 analytic and bounded for $k\in\,{\mathbb C}^{-}\setminus D$.
 \\
 \end{tabular}

 \vspace{2mm}
 An application of Cauchy's theorem yields
 \bea
 \int_{-\infty}^{\infty}{\rm
 e}^{ikx-\omega(k)t}\tilde{f}(t,k)dk&=&\int_{\partial D^{+}}{\rm
 e}^{ikx-\omega(k)t}\tilde{f}(t,k)dk,
 \nonumber \\
 \int_{-\infty}^{\infty}{\rm
 e}^{ik(x-L)-\omega(k)t}\tilde{g}(t,k)dk&=&-\int_{\partial D^{-}}{\rm
 e}^{ik(x-L)-\omega(k)t}\tilde{g}(t,k)dk,
 \nonumber
 \eea
 and (\ref{intrep1}) follows.
 
 We now derive the solution representation for example (a), (b) and 
 (c), with the boundary conditions considered in the previous section. 

 \subsubsection{Example (a)}
 For this example,   equation (\ref{intrep1}) becomes
\bea
q(x,t)&=&\frac{1}{2\pi}\left\{\int_{-\infty}^{\infty}{\rm
e}^{ikx-ik^{2}t}\hat{q}_{0}(k)dk+
 \int_{\partial D^{+}}{\rm e}^{ikx-ik^{2}t}(k\tilde{f}_{0}
 -i\tilde{f}_{1})dk\right.
 \nonumber \\
 &+&\left.
\int_{\partial D^{-}}{\rm e}^{ik(x-L)-
ik^{2}t}(k\tilde{g}_{0}
 -i\tilde{g}_{1})dk\right\},
\label{intrepa}
\eea
where $\partial D^{+}$ and $\partial D^{-}$ are the boundaries of the 
first and third 
quadrant of the complex $k$ plane, respectively. 

 \vspace{2mm}
  {\bf The Dirichlet problem}
  
  {\em 
  Let $q(x,t)$ satisfy equation (\ref{eq1}), the initial condition $q(x,0)=q_{0}(x)$ and 
  the boundary conditions $q(0,t)=f_{0}(t)$, $q(L,t)=g_{0}(t)$. 
  Then $q(x,t)$ admits the representation
  \bea
    q(x,t)&=&\frac{1}{2\pi}\left\{
     \int_{\partial D^{+}_{0}}{\rm
     e}^{ikx-ik^{2}t}\left[k\tilde{f}_{0}(t,k)+\frac{{\rm
     e}^{-ikL}N(k,t)-{\rm
     e}^{ikL}N(-k,t)}{\Delta(k)}\right]dk\right.
    \nonumber
    \\
    &&+\int_{\partial D^{-}_{0}}{\rm
     e}^{ik(x-L)-ik^{2}t}\left[k\tilde{g}_{0}(t,k)
     +\frac{N(-k,t)-N(k,t)}{\Delta(k)}\right]dk,
    \nonumber \\
    &&+
    \left.\int_{-\infty}^{\infty}{\rm
    e}^{ikx-ik^{2}t}\hat{q}_{0}(k)dk
       \right\},
       \label{eq2ex}
    \eea
  where $\hat{q}_{0}(k)$ is the Fourier transform of $q_{0}(x)$, 
 $N(k,t)$ is given by (\ref{Nk}), $\Delta(k)=\E-{\rm e}^{-ikL}$, and 
  $\partial D_{0}^{+}$ and $\partial D^{-}_{0}$ are the contours 
  $\partial D^{+}$ and $\partial D^{-}$ indented to pass above and below 
  the zeros of $\Delta(k)$ on the real axis, respectively.  
  }
  
  \vspace{2mm}
 Indeed, solving equations (\ref{A.2}) for $\tilde{f}_{1}$ and $\tilde{g}_{1}$ 
  and substituting the result in (\ref{intrepa}) we find
  \bea
  q(x,t)&=&\frac{1}{2\pi}\left\{
   \int_{\partial D^{+}}{\rm
   e}^{ikx-ik^{2}t}\left[k\tilde{f}_{0}(t,k)+\frac{{\rm
   e}^{-ikL}N(k,t)-{\rm
   e}^{ikL}N(-k,t)}{\Delta(k)}\right]dk\right.
  \nonumber
  \\
  &+&\int_{\partial D^{-}}{\rm
   e}^{ik(x-L)-ik^{2}t}\left[k\tilde{g}_{0}(t,k)
   +\frac{N(-k,t)-N(k,t)}{\Delta(k)}\right]dk,
  \nonumber \\
  &+&
    \int_{\partial D^{+}}{\rm
    e}^{ikx-ik^{2}t}\frac{{\rm
    e}^{-ikL}\hat{q}(t,k)-{\rm
    e}^{ikL}\hat{q}(t,-k)}{\Delta(k)}dk
    \label{eq2ex0}
   \\
   &+&\left.\int_{\partial D^{-}}{\rm
    e}^{ik(x-L)-ik^{2}t}\frac{\hat{q}(t,-k)-\hat{q}(t,k)}
    {\Delta(k)}dk+\int_{-\infty}^{\infty}{\rm
  e}^{ikx-ik^{2}t}\hat{q}_{0}(k)dk
     \right\},
   \nonumber
  \eea
  where $N(k,t)$ is given by (\ref{Nk}).
The determinant of (\ref{A.2}) is the function $\Delta(k)$ given in 
the statement. The 
zeros of this function are on the real line, which is part of 
the boundary  $\partial D$. Thus, when deforming the 
 contour to obtain the effective integral representation 
 (\ref{intrep1}), the
 part of the contour $\partial D^{+}$ along the positive real axis 
 must
 be deformed to a small circle above each of the points
 $k_{m}\in\,{\mathbb R}$, $m\in\,{\mathbb
 Z}^{+}$. Similarly the part
 $\partial D^{-}$ along the negative real axis must
 be deformed to a small circle below each of the points
 $k_{m}\in\,{\mathbb R}$, $m\in\,{\mathbb
 Z}^{-}$. It can then be 
 verified that the terms 
 involving the function $\hat{q}(t,k)$, after multiplication by ${\rm
e}^{-ik^{2}t}$, are analytic and bounded in 
 the indented domains
 $D^{+}$ and $D^{-}$ respectively. By an application of Jordan's 
 lemma, this implies that these terms give  a zero contribution ot 
 the representation.

 \subsubsection{Example (b)}
 For this example,   equation (\ref{intrep1}) becomes
\bea
q(x,t)&=&\frac{1}{2\pi}\left\{\int_{-\infty}^{\infty}{\rm
e}^{ikx-k^{2}t}\hat{q}_{0}(k)dk-
 \int_{\partial D^{+}}{\rm e}^{ikx-k^{2}t}(ik\tilde{f}_{0}
 +\tilde{f}_{1})dk\right.
 \nonumber \\
 &-&\left.
\int_{\partial D^{-}}{\rm e}^{ik(x-L)-
k^{2}t}(ik\tilde{g}_{0}
 +\tilde{g}_{1})dk\right\},
\label{intrepb}
\eea
where $\partial D^{+}$ and $\partial D^{-}$ are shown in figure \ref{heatdom}.

\vspace{2mm}
 {\bf The Dirichlet problem}

 {\em 
  Let $q(x,t)$ satisfy equation (\ref{eq1}), the initial condition $q(x,0)=q_{0}(x)$ and 
  boundary conditions $q(0,t)=f_{0}(t)$, $q(L,t)=g_{0}(t)$. 
  Then $q(x,t)$ admits the representation
  \bea
    q(x,t)&=&\frac{1}{2\pi}\left\{
     \int_{\partial D^{+}}{\rm
     e}^{ikx-k^{2}t}\left[k\tilde{f}_{0}(t,k)+\frac{{\rm
     e}^{-ikL}N(k,t)-{\rm
     e}^{ikL}N(-k,t)}{\Delta(k)}\right]dk\right.
    \nonumber
    \\
    &&+\int_{\partial D^{-}}{\rm
     e}^{ik(x-L)-k^{2}t}\left[k\tilde{g}_{0}(t,k)
     +\frac{N(-k,t)-N(k,t)}{\Delta(k)}\right]dk,
    \nonumber \\
    &&+
    \left.\int_{-\infty}^{\infty}{\rm
    e}^{ikx-k^{2}t}\hat{q}_{0}(k)dk
       \right\},
       \label{eq2ex2}
    \eea
  where $\hat{q}_{0}(k)$ is the Fourier transform of $q_{0}(x)$, 
 $N(k,t)$ is given by (\ref{Nk2}), and $\Delta(k)={\rm e}^{-ikL}-{\rm
e}^{ikL}$.
  }
  
  \vspace{2mm}
Indeed, solving
equations (\ref{agr2})-(\ref{sys2}) for $\tilde{f}_{1}$ 
and $\tilde{g}_{1}$, and
substituting the result in (\ref{intrepb}) we find (\ref{eq2ex2}). 

The solution of the system (\ref{agr2})-(\ref{sys2})
includes also terms involving the function $\hat{q}(t,k)$.
After multiplication by ${\rm
e}^{-k^{2}t}$, since 
$\Delta(k)$ has no zeros in $D$, these terms are analytic and 
bounded as
$k\to\infty$ in $D$.
An application of Jordan's
lemma implies that these terms give a zero contribution.

 \vspace{2mm}
 {\bf The Robin problem}
 
 {\em 
  Let $q(x,t)$ satisfy equation (\ref{eq2}), the initial condition $q(x,0)=q_{0}(x)$ and 
  boundary conditions (\ref{robin}). 
  Then $q(x,t)$ admits the representation
  $$
  q(x,t)=\frac{1}{2\pi}\left\{\int_{\partial D^{+}}{\rm
  e}^{ikx-k^{2}t}\left[(ik+\alpha)\frac{{\rm e}^{-ikL}(\beta+ik)N(-k,t)
  -{\rm e}^{ikL}(\beta-ik)N(k,t) 
  }{\Delta(k)}+\tilde{h_{1}}(t,k)
   \right]dk+
  \right.
  $$
  \be
  \left.
  \int_{\partial D^{-}}{\rm
   e}^{ik(x-L)-k^{2}t}\left[(ik+\beta)\frac{
  (\alpha+ik)N(-k,t)-(\alpha-ik)N(k,t)}{\Delta(k)}+\tilde{h_{2}}(t,k)
  \right]dk+\int_{-\infty}^{\infty}{\rm
  e}^{ikx-k^{2}t}\hat{q}_{0}(k)dk
  \right\},
  \label{eq4ex} \ee
  where $\hat{q}_{0}(k)$ is the Fourier transform of $q_{0}(x)$, 
 $N(k,t)$ is given by (\ref{Nheat}), $\Delta(k)$ is given by 
 (\ref{heatdet}) and 
 $$
 \tilde{h_{1}}(t,k)=\int_{0}^{t}{\rm e}^{k^{2}s}h_{1}(s)ds, \quad
 \tilde{h_{2}}(t,k)=\int_{0}^{t}{\rm e}^{k^{2}s}{\rm
 e}^{-ikL}h_{2}(s)ds.
 $$ 
  }
  
  \vspace{2mm}
 To derive equation (\ref{eq4ex}) we use the solutions
 (\ref{solveheat1}), (\ref{solveheat2}) of
 the system (\ref{grheat2})-(\ref{grheat3}).
 In the expressions (\ref{solveheat1}), (3.14) there appears
 also a term involving the function $\hat{q}(t,k)$. This term, when
 multiplied by ${\rm e}^{ikx-k^{2}t}$, is analytic and
  bounded for  $k\to\infty$ in $D^{+}$. Similarly, the term in 
 (\ref{solveheat2})
  involving the function $\hat{q}(t,k)$, when
 multiplied by ${\rm e}^{ik(x-L)-k^{2}t}$, is analytic and
  bounded for  $k\to\infty$ in $D^{-}$.
  Since the real axis is outside the domain $D$, $\Delta(k)\neq 0$ for 
  k in $D$, and an application of Jordan's
 lemma therefore  implies that these terms give a zero contribution.

\subsubsection{Example (c)}
In this case, the representation (\ref{intrep1}) becomes
\bea
q(x,t)&=&\frac{1}{2\pi}\left\{\int_{-\infty}^{\infty}{\rm
e}^{ikx-i(k-k^{3})t}\hat{q}_{0}(k)dk-
 \int_{\partial D^{+}}{\rm e}^{ikx-i(k-k^{3})t}((k^{2}-1)\tilde{f}_{0}
 -ik\tilde{f}_{1}-\tilde{f}_{2})dk\right.
 \nonumber\\
 &-&\left.
\int_{\partial D^{-}}{\rm 
e}^{ik(x-L)-(k^{2}-ik)t}((k^{2}-1)\tilde{g}_{0}
 -ik\tilde{g}_{1}-\tilde{g}_{2})dk\right\}.
\label{intrepc}
\eea
We now consider equation (\ref{eq3}) with the boundary conditions 
(\ref{bcc}).

\vspace{2mm}
{\em 
 Let $q(x,t)$ satisfy equation (\ref{eq3}), the initial condition $q(x,0)=q_{0}(x)$ and 
 boundary conditions (\ref{bcc}). 
 Then $q(x,t)$ admits the representation
 $$
 q(x,t)=\frac{1}{2\pi}\int_{\infty}^{\infty}{\rm
 e}^{ikx-i(k-k^{3})t}\hat{q}_{0}(k)dk+\frac{1}{2\pi}\int_{\partial 
 D^{+}}{\rm
 e}^{ikx-i(k-k^{3})t}\left((1-k^2)\tilde{f}_{0}(t,k)\right)dk
 $$
 \be
 +\frac{1}{2\pi}\int_{\partial D^{-}}{\rm
 e}^{ik(x-L)-i(k-k^{3})t}\left(1-k^2)\tilde{g}_0(t,k)+
 ik\tilde{g}_{1}(t,k)\right)dk
 \label{eq3ex}
 \ee
 $$
 +\frac{1}{2\pi}\int_{\partial D^{+}}{\rm
 e}^{ikx-i(k-k^{3})t}\left(ik\tilde{f}_{1}(t,k)+\tilde{f}_{2}(t,k)\right)dk
 +\frac{1}{2\pi}\int_{\partial D^{-}}{\rm
 e}^{ik(x-L)-i(k-k^{3})t}\tilde{g}_{2}(t,k)dk,
 $$
 where $\hat{q}_{0}(k)$ is the Fourier transform of $q_{0}(x)$, 
$\tilde{f}_{1}(k,t)$, $\tilde{f}_{2}(k,t)$, $\tilde{g}_{2}(k,t)$ are given 
by (1.16), $\partial D^{+}$ is the  branch of the hyperbola
$3{\rm Re}(k)^{2}-{\rm Im}(k)^{2}-1=0$  in the upper half plane, and $D^{-}$ is the branch of the 
same hyperbola 
in the lower half plane,
see figure \ref{fig1}.
  }
  
\vspace{2mm}
 Solving equations (\ref{agr3})-(\ref{sys3}) for $\tilde{f}_{1}$, $\tilde{f}_{2}$
 and $\tilde{g}_{2}$ and substituting the resulting expressions (see 
 equations (\ref{g23})) in
 (\ref{intrepc}),  we find (\ref{eq3ex}).

 We
 emphasise once more that the solution of the system (\ref{agr3})-(\ref{sys3})
 includes also the terms involving $\hat{q}(t,k)$. However, these terms do 
not  contribute to
 the solution, since 
 $\Delta(k)\neq 0$ for k in $D$, therefore an application of Jordan's
 lemma implies that their integral vanishes.

\subsubsection{Remarks on the equivalence of series and integral 
representations}
We consider example (b).
Equation (\ref{eq2ex}) shows that there exists a representation of
the solution of the Dirichlet problem for equation (\ref{eq2})
involving {\em only} integrals. It is of course possible to rewrite 
the
representation (\ref{eq2ex}) in the series form  using
Cauchy's theorem.  Indeed, it is easy to verify by studying the
boundedness of each exponential involved in (\ref{eq2ex})
that the
integrals along $\partial D^{+}$ and
$\partial D^{-}$ can be deformed to the real line, where the zeros
$k_{m}=\frac{m\pi}{L}$ of $\Delta(k)$ lie.
Computing explicitly the  residues at these poles, and manipulating
the result, it is easy to verify  that all integral terms cancel out, 
yielding 
\be
q(x,t)=
\frac{1}{2L}\sum_{m\in{\mathbb Z}}{\rm e}^{-k^{2}_{m}t}
 \left({\rm e}^{ik_{m}x}-{\rm e}^{-ik_{m}x}\right)
 \left(\hat{q}_{0}(k_{m})-\hat{q}_{0}(-k_{m})\right).
\label{seriesb}
\ee
The situation for example (a) is similar.

\section{Conclusions}\label{sec6}
We have illustrated the applicability of a {\em transform method} by
solving several concrete boundary value problems. It appears that the method is both
general and simple to implement. Indeed, the only mathematical tools
used in this paper are the Fourier transform and Cauchy's theorem.

An effort has been made to minimise technical considerations. In
particular, the given initial and boundary conditions are assumed to
be ``sufficiently smooth''. It is possible to work in a less
restrictive function class; for problems on the half line $x>0$
this is done in \cite{fsung2}, where general theorems are proven in
appropriate Sobolev spaces.

The general result about \bvps on an interval is that $q(x,t)$ can 
always be
expressed as an {\em integral in the complex $k$-plane}, see
equation (\ref{intrep1}). This integral
involves the Fourier transform $\hat{q}_{0}(k)$ of the initial
condition $q_{0}(x)$ and the $t$-transforms
$\{\tilde{f}_{j},\tilde{g}_{j}\}$ of {\em all} boundary values
$\{\partial_{x}^{j}q(0,t),\partial_{x}^{j}q(L,t)\}$. A subset of
these boundary values can be prescribed as boundary conditions. Thus
a subset $\tilde{f}_{p}$ and $\tilde{g}_{q}$, where $p$ takes $N$
values and $q$ takes $n-N$ values, can be computed immediately. The
remaining $\tilde{f}_{j}$'s and $\tilde{g}_{j}$'s can be expressed
through the solution of a system of $n$ algebraic equations obtained
from the global relation and from the equations derived from the
global relation by replacing $k$ with $\lambda_{j}(k),
\;j=1,\ldots,n-1$. The relevant expressions involve the unknown
function $\hat{q}(t,k)$ and the function $1/\Delta(k)$, where
$\Delta(k)$ is the determinant of the associated system of $n$
algebraic equations. However, using the integral representation
of $q(x,t)$ (equation (\ref{intrep1})) it can be shown that
{\bf (i)}:  If $\Delta(k)\neq 0$ in the domain $D$ (defined by
equation (\ref{Ddef})), then the contribution of the terms involving
$\hat{q}(t,k)$  vanishes.
{\bf (ii)}: If $\Delta(k)$ has zeros in $D^{+}$ (hence also in 
$D^{-}$), then the
contours $\partial D^{+}$ and $\partial D^{-}$ must be indented to
pass above or below these zeros. 
In these
particular cases, as well as in general for the case of even order PDEs, 
there exists an alternative representation
consisting only of an infinite series. The simplest way to obtain this 
representation is to evaluate (\ref{gr1n}) at the zeros of $\Delta(k)$. 

\vspace{2mm}
 The basic examples (\ref{eq1})-(\ref{eq3}) were chosen
in order to illustrate the above cases:

{\bf (a)} For equation (\ref{eq1}), the domain $D$ is the union of
the first and third quadrant of the complex $k$-plane. For
Dirichlet (or Neumann) boundary conditions,  $\Delta(k)=0$ if
$k=\frac{m\pi}{L}$, $m\in\,\mathbb Z$, 
thus the contours along the
positive and negative real axes must be indented. Furthermore, it
can be shown that the relevant integral can be computed entirely
in terms of a sum of residues. Thus there exists an alternative 
representation of the solution in the form of an infinite series. 
This representation is consistent with the well
known form of the solution as given by a sine series.

{\bf (b)} For equation ({\ref{eq2}), the domain $D$ is shown in
figure (\ref{heatdom}). For Dirichlet boundary conditions,
$\Delta(k)=0$ for $k=\frac{m\pi}{L}$, $m\in{\mathbb
Z}$, thus $\Delta(k)\neq 0$ for $k$ in $D$ and there is {\em no} need 
to indent the contour. However, since this equation is of even 
order, it is still possible to rewrite these integrals in terms of
an infinite sum, see equation (\ref{eq2ex}). This is consistent
with the classical series solution.

For  more complicated boundary conditions, such as Robin-type 
conditions, the integral representation of the solution (\ref{eq4ex}) is more 
{\em effective} than the series representation.

{\bf (c)} For equation ({\ref{eq3}), the domain $D$ is shown in
figure (\ref{fig1}). If the given  boundary conditions 
are
$q(0,t)$, $q(L,t)$ and $q_{x}(L,t)$, the determinant $\Delta(k)$ is
given by equation (\ref{detc1}). In can be shown that the zeros of
$\Delta(k)$ are asymptotically on the three lines $\{k:\;{\rm 
arg}(k)=\pi/6,
\;{\rm
arg}(k)=5\pi/6\;{\rm
arg}(k)=3\pi/2\}$, and that  $\Delta(k)$ has {\em no} zeros in $D$.
In this case the contour {\em cannot} be deformed to pick up 
the contribution of the residues at these zeros. Equivalently, the solution 
of the system obtained from the global relation {\em cannot} be evaluated at these zeros
to obtain directly a series representation for $q(x,t)$. This is the generic behaviour 
for odd order problems. These problems admit a series representation only for 
coupled boundary conditions (including the periodic case). 

\vspace{2mm} We note that even in the cases when it is possible to
express the integral representation of the solution in terms of an
infinite sum, the integral form may have some advantages. For
example, the integral representation, in contrast with the series
one, is {\em uniformly convergent}  both at $x=0$ and at $x=L$. In
addition, integrals are more convenient than sums for studying the
long time asymptotic behaviour of the solution. Furthermore, for
many concrete examples it is possible to use the integral
representation and Cauchy's theorem to compute $q(x,t)$
explicitly. It is interesting that in these computations, one does
{\em not} compute the relevant integral by using the residues
associated with $\Delta(k)=0$, but one precisely {\em avoids} these
zeros.

We emphasise that an additional advantage of the representation
(\ref{intrep1}) is that it does not require a detailed analysis of the
function $\Delta(k)$. This is to be contrasted with the classical
approach: since the zeros of $\Delta(k)$ define the discrete
spectrum of the associated $x$-differential operator,
a {\em detailed} characterisation of this set is crucial for the 
derivation
of the associated basis of eigenfunctions. This important
advantage of the representation (\ref{intrep1}) for the Robin problem 
for the heat equation is illustrated by comparinng equations (\ref{eq4ex}) 
and (\ref{hfin}).

\vspace{2mm}
We conclude with some remarks:

{\bf (1)}
Coupled boundary conditions can also be analysed using our
method.  For example, for the
\bvp\,
$$
q_{t}+q_{xxx}=0,\quad q(x,0)=q_{0}(x),\quad q(0,t)=q(L,t)=0,
q_{x}(L,t)+\alpha q_{x}(0,t)=0,
$$
the solution {\em can} be expressed in a series of eigenfunctions
\cite{zhang}. This is a consequence of the fact that the
associated $x$-differential operator with the above boundary
conditions, has a self-adjoint extension. In this case, it is possible
to deform the integrals appearing in our representation and to
rewrite them in terms of an infinite series. Alternatively, the series 
can be obtained by evaluting the relevant system derived from equations 
(\ref{gr1n}) at the zeros of $\Delta(k)$.

{\bf (2)} Earlier work on two-point boundary value problems for
linear evolution PDEs has appeared in \cite{fp4}, but the existence of
zeros of $\Delta(k)$ was not investigated. The
determination of conditions for well
posedness is presented in \cite{pel}. The role of the zeros of the 
determinant
$\Delta(k)$ and their determination in some particular cases are
studied in
\cite{pel2}.

\section*{Appendix}
\subsection*{A. Well posed problems - the determination of $N$}

In order to motivate
the choice of  $N$, we present two examples for each of 
equations (\ref{eq1}) and (\ref{eq3}).

\subsubsection{Equation (\ref{eq1})}

{\bf (a.1) ${\bf q(0,t)=f_{0}(t)}$, ${\bf q(L,t)=g_{0}(t)}$}

In this case, the two functions $\tilde{f}_{0}(t,k)$ and
$\tilde{g}_{0}(t,k)$ are known, and we view equations 
(\ref{agr1})-(\ref{sys1}) as a system for $\tilde{f}_{1}$ and $\tilde{g}_{1}$. We concentrate
only on the dependence of the solution on the term $\hat{q}(t,k)$.
Solving equations (\ref{A.2}) for $\tilde{f}_{1}$ and $\tilde{g}_{1}$ we find
$$
i\tilde{f}_{1}=(known)+{\rm e}^{ik^{2}t}\left[\frac{{\rm
e}^{2ikL}\hat{q}(t,k)-\hat{q}(t,-k)}{1-{\rm e}^{2ikL}}\right],
\nonumber \\
i\tilde{g}_{1}=(known)+{\rm
e}^{ik^{2}t}\left[\frac{\E\hat{q}(t,k)-\E\hat{q}(t,-k)}{1-{\rm
e}^{2ikL}}\right].
\eqno(A.1)
$$
The terms  involving the function $\hat{q}(t,k)$ in equation
(A.1(a)) are
{\em bounded in $k$} as $k\to\infty$
if $k\in\,D^{+}=D\cap {\mathbb C}^{+}$.

Indeed, for
$k$ in ${\mathbb C}^{+}$,
 $$
\frac{{\rm e}^{2ikL}\hat{q}(t,k)-\hat{q}(t,-k)}
{1-{\rm e}^{2ikL}}\sim {\rm e}^{2ikL}\hat{q}(t,k)-\hat{q}(t,-k),
$$
and both these terms are bounded in  ${\mathbb C}^{+}$. Since the term
${\rm e}^{ik^{2}t}$ is bounded in $D$, our claim follows.

Similarly, the terms involving the function $\hat{q}(t,k)$ in
equation (A.1(b)) are bounded as $k
\to\infty$ if $k\in\,D^{-}=D\cap {\mathbb C}^{-}$. Indeed,
 $$
\frac{\E\hat{q}(t,k)-\E\hat{q}(t,-k)}{1-{\rm
e}^{2ikL}}\sim {\rm e}^{-ikL}\hat{q}(t,k)- {\rm 
e}^{-ikL}\hat{q}(t,-k),
$$
and both terms are bounded in ${\mathbb C}^{-}$.

\vspace{2mm}
{\bf (a.2) ${\bf q(0,t)=f_{0}(t)}$, ${\bf q_{x}(0,t)=f_{1}(t)}$}

In this case we obtain a system for the two functions
$\tilde{g}_{0}(t,k)$ and
$\tilde{g}_{1}(t,k)$. Considering explicitly only the terms involving
the function $\hat{q}(t,k)$, we
find
\bea
k\tilde{g}_{0}&=&known+{\rm e}^{ik^{2}t}\left[\frac{-{\rm
e}^{2ikL}\hat{q}(t,k)+\hat{q}(t,-k)}{2\E}\right],
\nonumber \\
i\tilde{g}_{1}&=&(known)+{\rm
e}^{ik^{2}t}\left[\frac{{\rm
e}^{2ikL}\hat{q}(t,k)+\hat{q}(t,-k)}{2\E}\right].
\nonumber
\eea
In this case  {\em not all} the terms  containing the function
$\hat{q}(t,k)$ are bounded as $k\to\infty$ for
$k\in\,D^{-}=D\cap {\mathbb C}^{-}$.
For example,
$$
\frac{-{\rm e}^{2ikL}\hat{q}(t,k)+\hat{q}(t,-k)}{2{\rm
e}^{ikL}}=-\frac{\E\hat{q}(t,k)}{2}+{\rm e}^{-ikL}\hat{q}(t,-k)$$
and the second term is bounded   in ${\mathbb C}^{-}$, but the
first term is bounded  in ${\mathbb C}^{+}$. Thus in the above
example, $N=1$, i.e. one boundary condition must be
prescribed at each end of the interval.

\subsubsection{Equation (\ref{eq3})}
The first set of
boundary conditions, with $N=1$, yields a well posed problems.
The second, with $N=2$, does
not.

\vspace{2mm}
{\bf (c.1) ${\bf q(0,t)=f_{0}(t)}$, ${\bf q(L,t)=g_{0}(t)}$, ${\bf
q_{x}(L,t)=g_{1}(t)}$}

In this case the unknown functions
$\tilde{f}_{1}$, $\tilde{f}_{2}$ and $\tilde{g}_{2}$ are given by
\bea
&&i\tilde{f}_{1}=(known \;\;terms)
\nonumber \\
&&+{\rm e}^{i(k-k^{3})t}\left[\frac{
{\rm
e}^{-i\lambda_{1}L}(\hat{q}(k)-\hat{q}(\lambda_{2}))+{\rm
e}^{-i\lambda_{2}L}(\hat{q}(\lambda_{1})-\hat{q}(k))+
{\rm e}^{-ikL}(\hat{q}(\lambda_{1})-\hat{q}(\lambda_{2}))
}{\Delta(k)}\right]
\nonumber \\
&&\tilde{f}_{2}=(known \;\;terms)
\nonumber \\
&&+{\rm e}^{i(k-k^{3})t}\left[\frac{{\rm
e}^{-i\lambda_{1}L}(k\hat{q}(\lambda_{2})-\lambda_{2}\hat{q}(k))+{\rm
e}^{-i\lambda_{2}L}(\lambda_{1}\hat{q}(k)-k\hat{q}(\lambda_{1}))+
{\rm e}^{-ikL}
(\lambda_{1}\hat{q}(\lambda_{2})-\lambda_{2}\hat{q}(\lambda_{1}))}
{\Delta(k)}
\right]
\nonumber\\
&&\tilde{g}_{2}=(known \;\;terms)
\nonumber 
\eea
$$
+{\rm
e}^{i(k-k^{3})t}\left[\frac{\hat{q}(\lambda_{1})(k-\lambda_{1})+
\hat{q}(\lambda_{2})(\lambda_{2}-k)+
\hat{q}(k)(\lambda_{2}-\lambda_{1})}
{\Delta(k)}\right]
\eqno(A.2) 
$$
where
$\Delta(k)$ is given by equation (\ref{detc1}) and
to simplify the
notation we have suppressed the $t$ dependence in $\hat{q}(t,k)$.

The terms  containing the function $\hat{q}(t,k)$ in equations
(A.2(a),(b)) are
{\em bounded in $k$} as $k\to\infty$
if $k\in\,D^{+}=D\cap {\mathbb C}^{+}$. This follows from the
observation that if
$k\in D^{+}$, then $\lambda_{1}$ and $\lambda_{2}$ are in
${\mathbb C}^{-}$.
For example, consider the function $\tilde{f}_{1}(k)$. Since ${\rm e}^{ikL}$, 
${\rm e}^{-i\lambda_{1}L}$  and ${\rm e}^{-i\lambda_{2}L}$  are all bounded for $k$ in ${\mathbb C}^{+}$, 
the
bracket appearing on the right hand side of $i\tilde{f}_{1}$, as
$k\to\infty$, is asymptotically given by
$$
\frac{
({\rm e}^{i(k-\lambda_{1})L}-{\rm e}^{i(k-\lambda_{2})L})\hat{q}(k)
+(1-{\rm e}^{i(k-\lambda_{2})L})\hat{q}(\lambda_{1})
+({\rm
e}^{i(k-\lambda_{1})L}-1)\hat{q}(\lambda_{2})}{\lambda_{1}-\lambda_{2}}.
$$
All terms in this expression are bounded when $k\in\,{\mathbb C}^{+}$.
In addition, ${\rm
e}^{i(k-k^{3})t}$ is bounded for all $k
\in D$, and the claim follows.
Similarly, the terms containing $\hat{q}(k)$ in equation
(A.2(c)) are
{\em bounded in $k$} as $k\to\infty$
if $k\in\,D^{-}=D\cap {\mathbb C}^{-}$.

\vspace{2mm}
{\bf (c.2) ${\bf q(0,t)=f_{0}(t)}$, ${\bf q(L,t)=g_{0}(t)}$,
${\bf q_{x}(0,t)=f_{1}(t)}$}

In this case, the unknown functions $\tilde{f}_{2}$, $\tilde{g}_{1}$
and $\tilde{g}_{2}$ are given by
\bea
&&\tilde{f}_{2}=(known \;\;terms)
\nonumber \\
&&+{\rm e}^{i(k-k^{3})t}\left[\frac{{\rm
e}^{ikL}(\lambda_{2}-\lambda_{1})\hat{q}(k)+{\rm
e}^{i\lambda_{2}L}(\lambda_{1}-k)\hat{q}(\lambda_{2})+
{\rm
e}^{i\lambda_{1}L}(k-\lambda_{2})\hat{q}(\lambda_{1})}{\Delta(k)}\right]
\nonumber \\
&&i\tilde{g}_{1}=(known \;\;terms)
\nonumber \eea
$$
+{\rm e}^{i(k-k^{3})t}\left[\frac{{\rm
e}^{-i\lambda_{1}L}(\hat{q}(k)-\hat{q}(\lambda_{2}))+
{\rm
e}^{-i\lambda_{2}L}(\hat{q}(\lambda_{1})-\hat{q}(k))+
{\rm e}^{-ikL}(\hat{q}(\lambda_{2})-\hat{q}(\lambda_{1}))
}
{\Delta(k)}\right]
\eqno(A.3)
$$
\bea
&&\tilde{g}_{2}=(known \;\;terms)
 \nonumber\\
&&+{\rm
e}^{i(k-k^{3})t}\left[\frac{{\rm
e}^{-i\lambda_{1}L}\lambda_{1}(\hat{q}(\lambda_{2})-\hat{q}(k))
+{\rm
e}^{-i\lambda_{2}L}\lambda_{2}(\hat{q}(k)-\hat{q}(\lambda_{1}))
+{\rm e}^{-ikL}k(\hat{q}(\lambda_{1})-\hat{q}(\lambda_{2})
}
{\Delta(k)}\right]
\nonumber
\eea
where $\Delta(k)$ is given by
$$
\Delta(k)={\rm e}^{ikL}(\lambda_{2}-\lambda_{1})
+{\rm e}^{i\lambda_{2}L}(\lambda_{1}-k)+{\rm
e}^{i\lambda_{1}L}(k-\lambda_{2}).
$$
As for example (a.2),  {\em not all} the terms containing the
unknown function $\hat{q}(t,k)$ in equations (A.3)
are bounded for all $k\in\,D^{+}$ or $k\in\,D^{-}$. As an example,
consider the terms
in (A.3(b)), which  should be bounded as $k\to\infty$
for all  $k\in D^{-}$. Choose $k$ such that $\lambda(k)\in\,D^{-}$ and
$\lambda_{2}(k)\in\,D^{+}$. Then $\E$ and ${\rm e}^{i\lambda_{1}L}$ 
are
not bounded as $k\to\infty$, while ${\rm e}^{i\lambda_{2}L}$
is bounded. Since asymptotically $\lambda_1\sim {\rm e}^{2\pi i/3}k$, 
it is easy to verify
that, for $k\in \,D^-$ such that $\lambda_1(k)\in \,D^-$,
the dominant term in the denominator is ${\rm e}^{i\lambda_{1}L}$.
Hence (A.3(b)), as $k\to\infty$, is  given by
 $$
{\rm e}^{i(k-k^{3})t}
\left[{\rm
e}^{-2i\lambda_{1}L}(\hat{q}(k)-\hat{q}(\lambda_{2}))+
{\rm
e}^{ikL}(\hat{q}(\lambda_{1})-\hat{q}(k))+
{\rm e}^{i\lambda_2L}\hat{q}(\lambda_{2})-\hat{q}(\lambda_{1})\right]
 $$
and
the term $\E\hat{q}(k)\sim {\rm e}^{ik(L-x)}$
is {\em not} bounded for  $k\in\,D^{-}$
(to simplify the
notation we have again suppressed the $t$ dependence in 
$\hat{q}(t,k)$).

\subsection*{B. The zeros of finite exponential sums}

The zeros of the function (\ref{tildedet3}) coincide with the zeros 
of 
$$F(z)={\rm e}^{z}+\zeta{\rm e}^{\zeta z}+\zeta^{2}{\rm e}^{\zeta^{2} 
z}, \quad \zeta={\rm e}^{2\pi i/3},
$$
where $z=-ikL$, $k\in\mathbb C$.
Following the general theory given in \cite{levin}, one can use a simple 
geometric construction to characterise the distributions of the zeros of functions 
of the form
$$
G(z)={\rm e}^{z}+a_{1}{\rm e}^{\lambda _{1}z}+\ldots+a_{n}{\rm 
e}^{\lambda_{n}z},
\eqno(B.1)
$$
where $a_{i}$, $\lambda_{i}$ are complex constants, such that the 
$n$-polygon with vertices the points $1$, $\lambda_{1}$, \ldots, $\lambda_{n}$ is not 
degenerate. 
In this case, the zeros of this function are clustered along the rays 
emanating from the origin with direction orthogonal to the sides of 
the polygon, and 
can only accumulate at infinity on these rays (regardless of the values of the constants $a_{i}$).

For the function $F(z)$, the associated polygon is the triangle with vertices 
the third roots of unity, $1$, ${\rm e}^{2\pi i/3}$ and ${\rm e}^{4\pi 
i/3}$. The rays normal to the sides of this triangle have directions
$\pi/3$, $\pi$ and $-\pi/3$. In terms of the variable $k=iz/L$, these 
rays have directions $\pi/6$, $5\pi/6$ and $3\pi/2$. Therefore the 
zeros  of the function
(\ref{tildedet3}) are guaranteed to lie asymptotically along these 
rays. In this particular case, the added symmetry of the function 
implies that all zeros lie precisely on these rays (see \cite{pel2} 
for a direct proof). However, in general the knowledge of the 
asymptotic distribution of these zeros is sufficient for the present 
purposes.

We note that this argument can be applied to the problems of order $n$ 
by reducing the relevant determinant asymptotically to the simple form 
(B.1).

\subsection*{C. The classical transform approach}

The form of the particular solution $E(x,t)$ of the PDE suggests that the most 
convenient representation is the one obtained by a Fourier transform 
with respect to $x$.  However, as it was mentioned earlier, for 
boundary value problem for odd order equations with uncoupled boundary 
conditions, there does {\em not} exist an appropriate $x$-transform. 
For such problems, one can
use a Laplace transform with respect to $t$.
Indeed, the
particular solution $E(x,t)$ can be rewritten in the form
${\rm e}^{st+ik(s)x}$, where $k$ satisfies the $n$-th order equation
$$
s+\omega(k)=0.
$$
%
%
%
%

We show here how equation (\ref{eq3}) can be solved using this 
appproach.

Let $\tilde{q}(x,s)$ be the Laplace transform of $q(x,t)$, i.e. 
$$
\tilde{q}(x,s)=\int_0^{\infty} {\rm e}^{-st}q(x,t)dt,\quad {\rm 
Re}\,(s)>0.
\quad \eqno(C.1)
$$
Applying  the Laplace transform to equation (\ref{eq3}) we find
$$
\tilde{q}_{xxx}+\tilde{q}_x+s\tilde{q}=q_0(x), \quad {\rm Re}\,(s)>0.
\quad \eqno(C.2)
$$
The solutions of the homogeneous version of this equation are given by
$$
\tilde{q}(x,s)={\rm e}^{\lambda x}, \quad \lambda^3+\lambda+s=0.
\quad \eqno(C.3)
$$
 We distinguish the three roots of the cubic equation (C.3) by their 
large $s$ behaviour:
$$
\lambda_1\sim -s^{\frac{1}{3}},\quad \lambda_2\sim -\alpha 
s^{\frac{1}{3}},
\quad \lambda_3\sim -\alpha^2s^{\frac{1}{3}}, \quad \alpha={\rm 
e}^{2\pi i/3}.
\quad \eqno(C.4)
$$
If $-\pi/2<{\rm arg}(s)<\pi/3$ then, for large $s$,
$$
-\frac{\pi}{6}<{\rm arg}(s)<\frac{\pi}{6},\quad 
\frac{\pi}{2}<{\rm arg}(\alpha s)<\frac{5\pi}{6},\quad
\frac{7\pi}{6}<{\rm arg}(\alpha^2s)<\frac{3\pi}{2}.
$$
Thus, for large $s$,
$$
{\rm Re}\,(\lambda_1)<0,\quad {\rm Re}\,(\lambda_j)>0,\;\;j=2,3.
\quad \eqno(C.5)
$$
A solution of the inhomogeneous equation (C.2) is given by
$$
\tilde{q}(x,s)=\sum_{j-1}^3 \beta_j\int_0^x {\rm 
e}^{\lambda_j(x-\xi)}q_0(\xi)d\xi,
\quad \eqno(C.6)
$$
where the constants $\{\beta_j\}_1^3$ satisfy the algebraic conditions
$$
\sum_1^3 \beta_j=0,\quad
\sum_1^3 \lambda_j\beta_j=0,\quad\sum_1^3(1+\lambda_j^2) \beta_j=0.
\quad \eqno(C.7)
$$
In the integral 
$$\int_0^x {\rm e}^{\lambda_j(x-\xi)}q_0(\xi)d\xi
$$
we have  $x-\xi>0$, thus the integrand is bounded as $s\to\infty$, 
if Re$(\lambda_j)<0$. On the other hand, 
if we replace $\int_0^x$ with $-\int_x^{\infty}$, then $x-\xi<0$, and 
the integrand  is bounded if Re$(\lambda_j)>0$. Hence, equations 
(C.5) imply that we must choose the following solution of equation 
(C.2):
$$
\tilde{q}(x,s)=c_1{\rm e}^{\lambda_1(x-L)} +c_2{\rm e}^{\lambda_2x} 
+c_3{\rm e}^{\lambda_3x}
$$
$$
-\beta_1\int_x^L {\rm e}^{\lambda_1(x-\xi)}q_0(\xi)d\xi
+\beta_2\int_0^x{\rm e}^{\lambda_2(x-\xi)}q_0(\xi)d\xi+
\beta_3\int_0^x{\rm e}^{\lambda_3(x-\xi)}q_0(\xi)d\xi,
\quad \eqno(C.8)
$$
where the $c_j$'s are constants independent of $x$.

Let $\tilde{f}_0$, $\tilde{g}_0$ and $\tilde{g}_1$ denote the Laplace 
transforms of the given boundary conditions $f_0(t)$, $g_0(t)$ and 
$g_1(t)$ respectively,
$$
\tilde{f}_0(s)=\int_0^{\infty}{\rm e}^{-st}f_0(t)dt,\quad 
\tilde{g}_0(s)=\int_0^{\infty}{\rm e}^{-st}g_0(t)dt,\quad 
\tilde{g}_1(s)=\int_0^{\infty}{\rm e}^{-st}g_1(t)dt. 
\quad \eqno(C.9)
$$
The definitions (C.1) and (C.9) imply
$$
\tilde{q}(0,s)=\tilde{f}_0(s), \quad
\tilde{q}(L,s)=\tilde{g}_0(s), \quad
\tilde{q}_x(L,s)=\tilde{g}_1(s).
\quad \eqno(C.10)
$$
Using equation (C.8) to evaluate $\tilde{q}(0,s)$, 
$\tilde{q}(L,s)$ and $\tilde{q}_x(L,s)$ we find the following set of 
three algebraic equations for $\{c_j\}_1^3$:
$$
{\rm 
e}^{-\lambda_1L}c_1+c_2+c_3=\beta_1\hat{q}_0(-\lambda_1)+\tilde{f}_0(s),
$$
$$
c_1+{\rm e}^{\lambda_2L}c_2+{\rm e}^{\lambda_3L}c_3=-\beta_2{\rm 
e}^{\lambda_2L}\hat{q}_0(-\lambda_2)-\beta_3{\rm 
e}^{\lambda_3L}\hat{q}_0(-\lambda_3)
+\tilde{g}_0(s),
$$
$$
\lambda_1c_1+\lambda_2{\rm e}^{\lambda_2L}c_2+\lambda_3
{\rm e}^{\lambda_3L}c_3=
-\beta_2\lambda_2{\rm e}^{\lambda_2L}\hat{q}_0(-\lambda_2)-
\beta_3\lambda_3{\rm e}^{\lambda_3L}\hat{q}_0(-\lambda_3)
+\tilde{g}_1(s),
\quad \eqno(C.11)
$$
where 
$$
\hat{q}_0(\lambda_j)=\int_0^L{\rm e}^{\lambda_j\xi}q_0(\xi)d\xi.
\quad \eqno(C.12)
$$
The determinant of the system (C.11) is given by
$$
\Delta=(\lambda_3-\lambda_2){\rm e}^{(\lambda_2+\lambda_3-\lambda_1)L}
+(\lambda_1-\lambda_3){\rm e}^{\lambda_3L}+
(\lambda_2-\lambda_1){\rm e}^{\lambda_2L}.
\quad \eqno(C.13)
$$
It can be shown that, as $s\to\infty$, the zeros of $\Delta$ are on 
the negative real axis. Indeed, recall that the zeros of
the analogous determinant in the complex $k$ plane   lie,  for large 
$k$, 
on the half lines arg$(k)=\pi/6$, arg$(k)=5\pi/6$ and 
arg$(k)=3\pi/2$.
Letting
$
s=i(k^3-k)$, the determinant (\ref{detc1}) reduces to the determinant 
(C.13), and furthermore the above three rays are mapped to the 
negative real 
axis in the complex $s$ plane. 
Solving equations (C.11) for the $c_j$'s and replacing the 
resulting expressions in equation (C.8) we find an expression for 
$\tilde{q}(x,s)$ which, as $s\to\infty$, can be shown to be bounded 
for Re$(s)>0$, Im$(s)\neq 0$. Actually, the relevant proof is 
identical to the proof of section \ref{sec4} which establishes that 
for this problem well-posedness demands one boundary condition at 
$x=0$ and two boundary conditions at $x=L$. 

\vspace{2mm}
    {\bf Remark C.1} The application of the Laplace transform involves the solution of 
two sets of algebraic equations, namely equations (C.7) and (C.11), 
while the application of the transform method used in this paper uses 
the solution of only one set of algebraic equations. Furthermore, the 
complex representation is based on the three roots of the {\em 
cubic} equation (C.3), while the representation (\ref{intrep1}) is 
based on the two roots of the {\em quadratic} equation 
$$\lambda^2+\lambda k+k^2-1=0.$$}

\vspace{2mm}
    {\bf Remark C.2} In order to rigorously justify the inversion 
    formula for the Laplace transform, usually the given boundary data 
    are not allowed to grow faster than linearly in the exponential. 
    On the other hand, the method presented in this paper
    does not require such restriction. 
  
\section*{Acknowlegdements}
The authors wish to thank A.R. Its, D.J. Needham and D. Powers for many useful
discussions.

BP gratefully acknowledges the support of the Nuffield foundation 
grant
NAL/00608/GR.



\begin{thebibliography}{1}


\bibitem{akhie}
N.I.Akhiezer and I.M.Glazman,
\newblock
{\em Theory of linear operators in Hilbert spaces},
\newblock
Nauka, 1966.

\bibitem{cg}
T. Colin and J.M. Ghidaglia,
\newblock
An initial-boundary value problem for the Korteweg-deVries equation
posed on a finite interval,
\newblock {\em Adv. in Diff. Eq.}, {\bf 6(12)} 2001, 1463-1492.

waves,

\bibitem{dym}
H. Dym and  H.P. McKean,
\newblock
{\em Fourier series and integrals},
\newblock Academic Press, 1972.

\bibitem{fok}
A.S. Fokas,
\newblock A unified transform method for solving linear and certain
nonlinear
  {PDE}'s,
\newblock {\em Proc. Royal Soc. Series A}, {\bf 453} 1997, 1411--1443.

\bibitem{best}
A.S. Fokas,
\newblock Two dimensional linear {P}{D}{E}'s in a convex polygon
\newblock {\em Proc. Royal Soc. Series A}, {\bf  457} 2001, 371--393.

\bibitem{ima1}
A.S. Fokas,
\newblock A new transform method for evolution PDEs,
\newblock {\em IMA J. Appl. Math}, {\bf 67} 2002, 559-590.

%
%
%

\bibitem{fp4}
A.S. Fokas and B.~Pelloni,
\newblock Two-point boundary value problems for linear evolution
equations
\newblock {\em Proc. Camb. Phil. Soc.} {\bf 17} 2001,
919-935.


\bibitem{fsung2}
A.S. Fokas and L.Y. Sung,
\newblock Initial boundary value problems for linear evolution
equations on the
  half line,
\newblock  {\em Ann. Math.} (in press).

\bibitem{levin}
B. Ja. Levin,
\newblock {\em Distribution of zeros of entire functions},
\newblock Translation of Mathematical Monographs, AMS, 1972.

\bibitem{pel}
B.~Pelloni,
\newblock Well-posed boundary value problems for
linear evolution equations on a finite interval,
\newblock  Math. Proc. Camb. Phil. Soc. {\bf 136} 2004,
361-382.

\bibitem{pel2}
B.~Pelloni.
\newblock
The spectral representation of two-point boundary value problems for
linear PDEs
\newblock submitted to Proc. Roy. Soc. Lond. A (2004)

\bibitem{tit}
C.E. Titchmarsh,
\newblock {\em The theory of functions},
\newblock Oxford University Press, 1962.

\bibitem{zhang}
B.Y. Zhang,
\newblock  Exact boundary controllability of the Korteweg-de
   Vries equation,
\newblock {\em SIAM J. Control Optim.} {\bf 37} 1999, 543-565 .

\end{thebibliography}
\end{document}